\documentclass[a4paper,12pt]{amsart}
\usepackage{amsfonts,amsthm,amssymb}
\usepackage{epsfig}
\usepackage{pstricks,slashbox,multirow}

\theoremstyle{plain}
\newtheorem{theorem}{Theorem}[section]

\newtheorem{definition}[theorem]{Definition}

\begin{document}

\date{}

\title[4CT - A New Simple Proof by Induction]{The Four Color Theorem - A New Simple Proof by Induction}

\author{V. Vilfred Kamalappan} 

\address{Department of Mathematics, Central University of Kerala 
\\ Periye - 671 320, Kasaragod, Kerala, India}

\email{vilfredkamalv@cukerala.ac.in}

\dedicatory{Dedicated to the memory of our parents, \\ Kamalappan Thomas and Mary Kamalappan  \linebreak and son Lerroy Wilson}


\begin{abstract}
In 1976, Appel and Haken achieved a major break through by proving the four color theorem $(4CT)$. Their proof is based on studying a large number of cases for which a computer-assisted search for hours is required. In 1997, Robertson, Sanders, Seymour and Thomas reproved the 4CT with less need for computer verification. In this paper, we present a simple proof to the $4CT$ based on mathematical induction and contraction. We consider, in our proof, possible colorings of a minimum degree vertex, its adjacent vertices and adjacent vertices of these adjacent vertices of simple planar graphs. 
\end{abstract}

\maketitle

\section{Introduction}

The Four Color Theorem (4CT) plays a prominent role in Graph Theory. The 4CT states that the vertices of every planar graph can be colored with at most four colors in such a way that any two adjacent vertices have different colors \cite{g08}, \cite{h90}-\cite{rs}. Throughout the paper, coloring of vertices means proper coloring of vertices so that adjacent vertices have different colors. 

A \textit{coloring} of a graph is an assignment of one color to each vertex so that no two adjacent vertices have the same color. The set of all vertices with any one color is independent and is called a \textit{color class}. An \textit{$n$-coloring} of a graph $G$ uses at most $n$ colors; it thereby partitions \textit{$V(G)$}, the set of vertices of $G$, into different color classes.

Heawood \cite{h90} showed that every planar graph is five colourable, Grotzsch \cite{g58} proved that any planar graph without triangle is three colorable and Grunbaum \cite{g63} proved that a planar graph with no more than three 3-circuits is three colorable.

In 1976, Appel and Haken \cite{ah1, ah2} achieved a major break through by proving the 4CT. Their proof is based on studying a large number of cases for which a computer-assisted search for hours is required. In 1997, Robertson, Sanders, Seymour and Thomas \cite{js} reproved the 4CT with less need for computer verification.

In this paper, we give a simple proof to the 4CT based on mathematical induction and contraction. An outline of the proof is given below and a detailed proof is presented in the next section. {\color{blue} We reached this stage of proof by making corrections/modifications and improvements several times based on reports of different referees at different times \cite{v87, v23}.} My sincere thanks to each one of these referees.  

For all basic ideas in graph theory, we follow \cite{h69}. The following definition is used in the proof.

\begin{definition} Let $G$ be a graph, $k\in\mathbb{N}$ and $v\in V(G)$. Then, the {\em $k$ distance neighbourhood of $v$} in $G$ is given by $N_k(v) = \{u\in V(G)/ 1 \leq d(u, v) \leq k \}$, the {\em $k$ distance closed neighbourhood of $v$} is $N_k[v]$ = $N_k(v) \cup \{ v \}$ and the set of $k^{th}$ distance vertices of $v$ is denoted by $\partial N_k(v)$ and is given by $\partial N_k(v)$ = $\{u\in V(G)/ d(u, v) = k \}$. Clearly, $N_1(v)$ = $N(v)$ and $N_1[v]$ = $N[v]$. 
\end{definition}

\vspace{.2cm}
\noindent
{\bf Outline of the proof.} 
\begin{enumerate}
\item [\rm (i)] Principle of Mathematical Induction is applied on $n$, the order of simple planar graphs to prove the theorem, $n\in\mathbb{N}$.

\item [\rm (ii)]  Assume that the theorem is true for simple planar graphs of order less than or equal to $n$. 

\item [\rm (iii)] Our aim is to prove the result on $H$, an arbitrary simple planar graph of order $n+1$.

\item [\rm (iv)] Consider a plane triangulation graph (maximal planar graph) $G$ of order $n+1$ of $H$ and prove the result on $G$ and thereby on $H$. 

\item [\rm (x)] {\color{blue} We classify graphs $G$ into two, one with a 3-cycle, say $(u_1, u_2, u_3)$ which is not a face boundary  and the other case, graphs $G$ without 3-cycle other than face boundaries.

 In the first case, we prove the result easily by considering proper coloring of proper subgraphs $G_O$ and $G_I$ of $G$ with interior and exterior of $(u_1, u_2, u_3)$, respectively and each with $(u_1, u_2, u_3)$}. 

 {\color{blue} In the second case, we consider the following.

\item [\rm (a)] $G$ contains a minimum degree vertex, say $v$, such that $d(v) \leq 5$.} {\color{red} By our assuption graph $G-v$ is 4 colorable. We consider all possible proper colorings of $G-v$ with 4 colors to prove the result. Here, coloring of vertices is done in such a way that vertices which are far away from $v$ is colored at first and then the next set of distant vertices from $v$ and so on.}

\item [\rm (b)] Let $n \geq 6$ and $v_1,v_2,\ldots,v_{d(v)}$ be the vertices adjacent to $v$ in $G$. They lie on a cycle, say $(v_1, v_2, . . . , v_{d(v)})$ in $G$ since $G$ is a maximal planar graph. 

\item [\rm (c)] Let $G_1$ = {\color{blue}  $G - N[v]$.} That is $G_1$ $=$ $G - \{v,v_1,v_2,\ldots,v_{d(v)}\}$.

\item [\rm (d)] {\color{blue} For a proper coloring of $G - v$, let $R_i$ be the set of colors of vertices of $G_1$ adjacent to $v_i$, $1 \leq i \leq d(v)$.} {\color{red} Then, for $i \neq j$, we say that sets $R_i$ and $R_j$ are adjacent whenever $v_i$ and $v_j$ are adjacent in $G-v$, $1 \leq i,j \leq d(v)$. In $G-v$, color of $v_i$ depends on the elements of $R_i$ and the colors of $v_{i-1}$ and $v_{i+1}$, $1 \leq i-1,i,i+1 \leq d(v)$.  }

{\color{blue} Thus, corresponding to a particular proper coloring of $G - v$ with four colors, we get a configuration (figure or coloring) consists of a set of values of $R_1, R_2, . . . , R_{d(v)}$ and the vertices $v_1, v_2, . . . , v_{d(v)}$ with their colors.} We consider all {\color{red} possible configurations corresponding to different possible proper colorings of $G-v$ with four colors to  prove the theorem.} 

In a {\color{blue} configuration, corresponding to a proper coloring of $G-v$ with four colors}, if cycle $(v_1, v_2, . . . , v_{d(v)})$ takes all the 4 colors of $G-v$, then with the same set of values of $R_1, R_2, \ldots,$ $R_{d(v)}$ and by making contraction on {\color{red} two edges of the cycle $(v_1, v_2, . . . , v_{d(v)})$ with a special property, then the resultant contracted graph, say $G_c$ of $G-v$, has to take more than four colors}. Thereby, the coloring is rejected and by assumption there exists an alternate coloring to $G_c$; thereby an alternate set of values to $R_1,$ $R_2,$ $\ldots,$ $R_{d(v)}$ in $G_c$ and thereby in $G_1$ (and in $G-v$). With these alternate values, we show, vertices $v_1, v_2, . . . , v_{d(v)}$ take less than 4 colors out of the 4 colors of $G-v$. 

Thus, {\color{red} we show that} there exists {\color{blue} proper} coloring to $G-v$ such that $G-v$ is 4 colorable while its vertices $v_1, v_2, ... , v_{d(v)}$ take less than 4 colors out of the 4 colors. And hence $G$ is four colorable. Thereby $H$ is four colorable.  \hfill $\Box$
\end{enumerate}

\section{Main Result}         
In this section, we present our proof to the four color theorem. Throughout the proof, graphs are considered in their planar form only.

\begin{theorem} [The Four Color Theorem (4CT)]\quad The vertices of every planar graph can be colored with at most four colors in such a way that any two adjacent vertices have different colors. 
\end{theorem}
\begin{proof}\quad It is enough to prove the theorem for any simple planar graph. The proof is based on the Principle of Mathematical Induction on the order $n$ of \textit{simple planar graphs}, $n\in\mathbb{N}$. 

The theorem is true for $n$ = 1, 2, 3, 4, 5. 

Assume the theorem for simple planar graphs of order less than or equal to $n$. Our aim is to prove the result for any simple planar graph of order $n+1$.

Let $H$ be an arbitrary simple planar graph of order $n+1$ and $G$ be a plane triangulation (maximal planar graph) of $H$. Now, it is enough to prove the result for $G$ since $H$ is a subgraph of $G$. 

\vspace{.2cm}
\noindent
{\it Claim}. Maximal planar graph $G$ of order $n+1$ is four colorable.

\vspace{.2cm}
{\color{blue} The following two cases of graph $G$ arise.}

\vspace{.2cm}
\noindent
{\bf {Case 1:}}\quad $G$ contains a 3-cycle that is not a face boundary.


{\color{blue} Let $(u_1, u_2, u_3)$ be a 3-cycle that is not a face boundary contained in $G$. Let $G_I$ be the proper subgraph of $G$ consisting of the 3-cycle and its interior and $G_O$ be the proper subgraph consisting of the 3-cycle and its exterior.} Thus, cycle $(u_1, u_2, u_3)$ is common to both $G_I$ and $G_O,$ $G_I \cup G_O$ = $G$ and $G_I \cap G_O$ $=$ $(u_1, u_2, u_3)$. {\color{blue} See Figures 1.1, 1.2 and 1.3. Proper subgraphs $G_I$ and $G_O$ of $G$, given in Figure 1.3, are presented separately in Figure 1.3a.}

\begin{figure}[ht]
    \centerline{\includegraphics[width=5in]{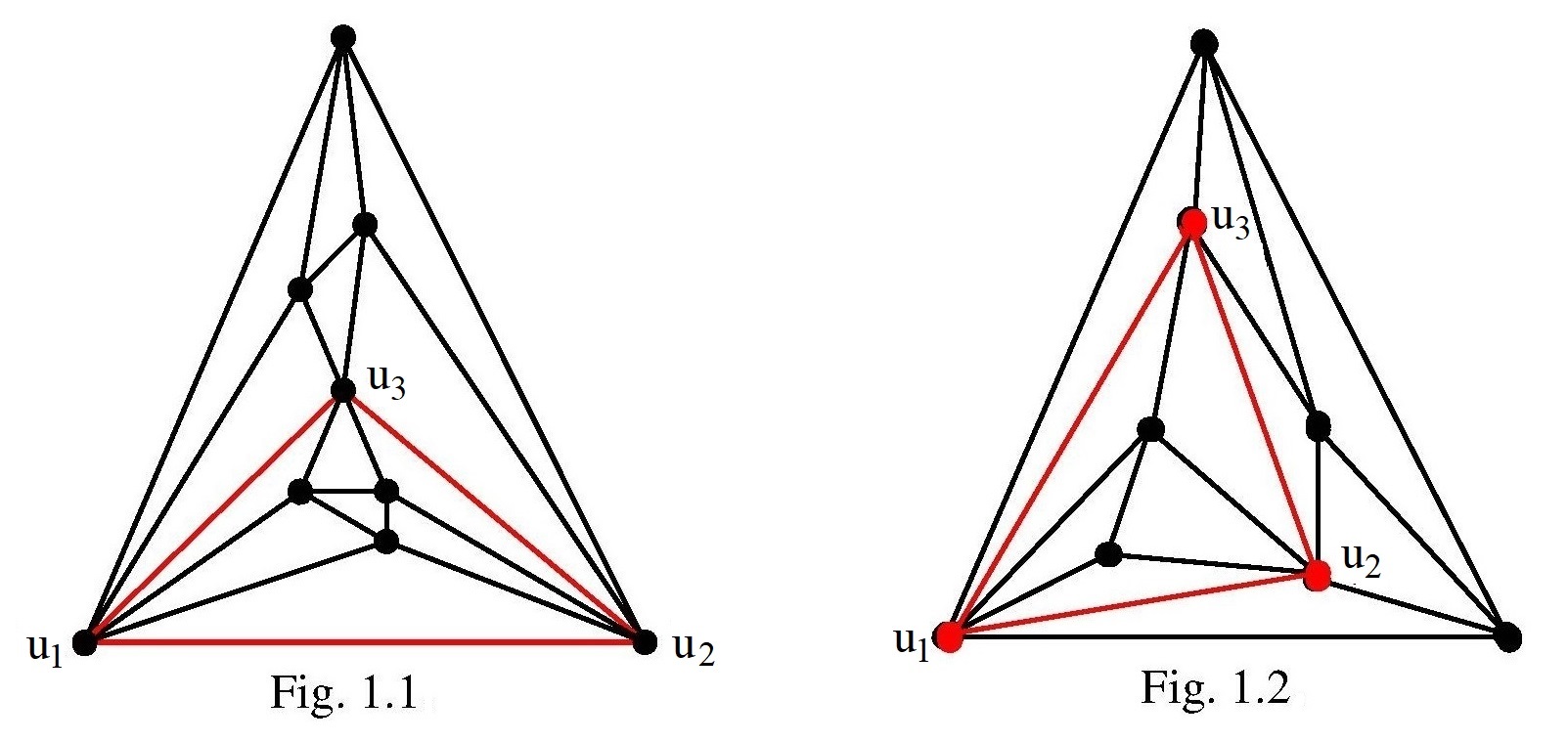}}
\end{figure}

\vspace{.5cm}
\begin{figure}[ht]
	\centerline{\includegraphics[width=5in]{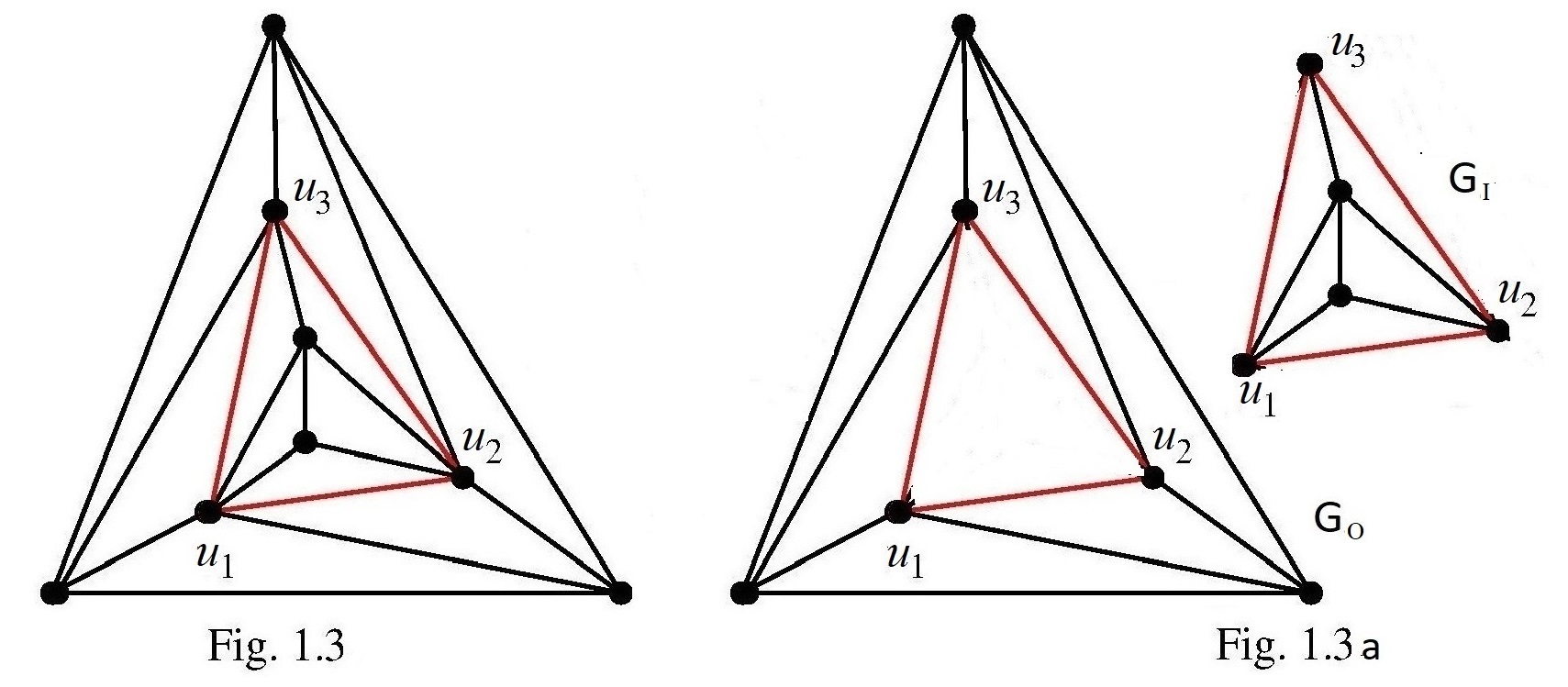}}
\end{figure}

Using the assumption, graphs $G_I$ and $G_O$ are four colorable, separately. Now, adjust the colors of vertices of $G_O$ by changing its color classes so that the colors of vertices $u_1,$ $u_2$ and $u_3$ of $G_O$ are same as in $G_I$ and thereby the combined graph $G$ of $G_I$ and $G_O$ is four colorable. For example, if $u_1,$ $u_2$ and $u_3$ are assigned with colors 1, 2, 3 in $G_I$ and 3, 4, 1 in $G_O$, respectively, then change the colors of vertices of $G_O$ by $3 \rightarrow 1$ (previous color of color class 3 of $G_O$ is replaced by color of color class 1 of $G_I$), $4 \rightarrow 2,$ $1 \rightarrow 3$ and $2 \rightarrow 4$. Now, $G$ is four colorable. 

\vspace{.2cm}
\noindent
{\bf Case 2:}\quad $G$ doesn't contain {\color{blue} 3-cycle other than face boundaries.}

Let $v$ be a vertex of $G$ with minimum degree. This implies, $d(v) \leq 5$ \cite{h69}. If $d(v) \leq 3$, then by the induction hypothesis $G-v$ is four colorable and since $v$ is adjacent to at the most 3 vertices of $G-v$, $v$ can be colored by a color among the four colors, different from the colors of its adjacent vertices. Hence, the theorem is true when $d(v) \leq 3$. When $d(v)$ = 4 or 5, {\color{blue} we consider the following to prove our claim.}

Let $v_1,v_2,\ldots,v_{d(v)}$ be the vertices adjacent to $v$ in $G$, $d(v)$ = $4$ or $5.$ In $G$, the vertices $v_1,v_2,\ldots,v_{d(v)}$ lie on a cycle, say $(v_1, v_2, . . . , v_{d(v)})$ since $G$ is a maximal planar graph. 
Figures 2.1 and 2.2 show the structure of graph $G$ as examples but not the whole graph.

\begin{figure}[ht]
	\centerline{\includegraphics[width=5in]{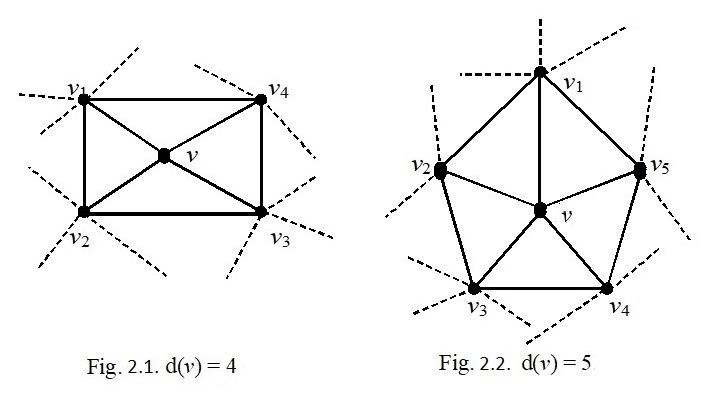}}
\end{figure}

Let $G_1$ = $G - N[v]$. That is $G_1$ = $G - \{v,v_1,v_2,\ldots,v_{d(v)}\}$.

{\color{blue} For a proper coloring of $G - v$ with four colors, let $R_i$ be the set of colors of vertices of $G_1$ adjacent to $v_i$, $1 \leq i \leq d(v)$. Then, for $i \neq j$, $R_i$ and $R_j$ are said to be {\textit adjacent} if and only if $v_i$ and $v_j$ are adjacent {\color{red} in $G-v$}, $1 \leq i,j \leq d(v)$. Figures 3.1 and 3.2 show $R_i$s and $v_i$s of $G-v$ for $d(v) = 4$ and $d(v)$ = 5, respectively. See Figures 3.1 and 3.2. In $G-v$, color of $v_i$ depends on the elements of $R_i$ and the colors of $v_{i-1}$ and $v_{i+1}$, $1 \leq i-1,i,i+1 \leq d(v)$.}

\begin{figure}[ht]
	\centerline{\includegraphics[width=5.5in]{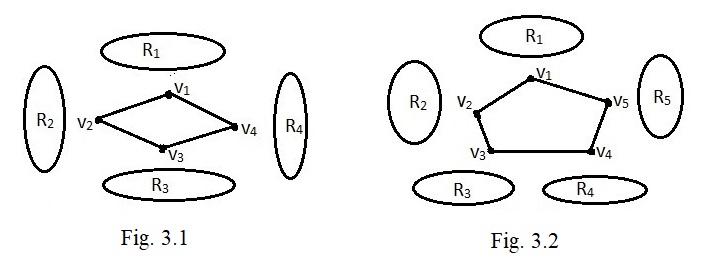}}
\end{figure}

{\color{red} In this proof, by our assuption graph $G-v$ is 4 colorable and to prove the result, we consider all possible proper colorings of $G-v$ with 4 colors. Here, coloring of vertices is done in such a way that vertices which are far away distance from $v$ is colored at first and then the next set of distant vertices from $v$ are colored and so on.}  
		
{\color{blue} Thus, corresponding to a particular proper coloring of $G - v$ with four colors, we get a configuration (figure or coloring) consists of a set of values of $R_1, R_2, . . . , R_{d(v)}$ and the vertices $v_1, v_2, . . . , v_{d(v)}$ with their colors.} We consider all {\color{red} possible configurations corresponding to different possible proper colorings of $G-v$ with four colors to  prove the theorem.} 
	
The following steps will reduce the number of possible configurations to Figures 7.1 to 7.5 when $d(v) = 4$ and to Figures 9.1 to 9.88 when $d(v) = 5$. 

\vspace{0.1in}
\noindent
{\bf Description of Figures 7.1 to 7.5 and 9.1 to 9.88}

{\color{blue} We consider proper colorings of $G-v$ with four colors.}

\begin{enumerate}
\item [{\rm 1.}] {\color{blue} In any proper coloring of $G-v$ with four colors}, $R_1$ $\cup$ $R_2$ $\cup$ $\ldots$ $\cup$ $R_{d(v)}$ $\subseteq$ $\{1, 2, 3, 4\}$ by assumption. 

\item [{\rm 2.}] {\color{red} $R_1$ $\cup$ $R_2$ $\cup$ $\ldots$ $\cup$ $R_{d(v)}$ = $\{1, 2, 3, 4\}$ can be avoided in $G-v$ by the following.
	 
  Suppose $R_1$ $\cup$ $R_2$ $\cup$ $\ldots$ $\cup$ $R_{d(v)}$ = $\{1, 2, 3, 4\}$ in $G_1$ for a particular proper coloring of $G-v$ with four colors. Now, in $G-v$, contract $v_1, v_2, \ldots, v_{d(v)}$ to a single vertex, say $v_c$ and let the contracted graph be $G_2$. Then, $v_c$ requires $5^{th}$ color in the contracted graph $G_2$ whose order is less than $n$. This is a contradiction to our assumption. This implies, there always exists proper coloring to $G-v$ with four colors  such that $R_1$ $\cup$ $R_2$ $\cup$ $\ldots$ $\cup$ $R_{d(v)}$ $\subseteq$ $\{1, 2, 3\}$.} Hence, we consider $R_1$ $\cup$ $R_2$ $\cup$ $\ldots$ $\cup$ $R_{d(v)}$ $\subseteq$ $\{1, 2, 3\}$.

 {\color{blue} Figures 4.1 to 4.4 are presented to get an idea of the structure of graphs $G-v$ and $G_1$ with $R_i$s and $v_i$s {\color{red} with $d(v)$ = 5 and $1 \leq i \leq d(v)$}. See Figures 4.1 to 4.4. Figure 4.5 shows the contracted graph $G_2$ of $G-v$ corresponding to Fig. 4.4. See Figure 4.5.} 
 \begin{figure}[ht]
 	\centerline{\includegraphics[width=5.5in]{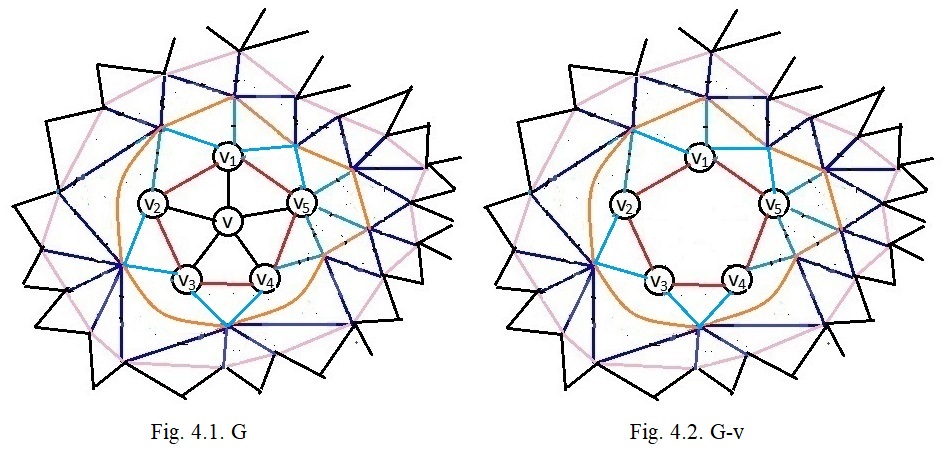}}
 \end{figure}
 \begin{figure}[ht]
 	\centerline{\includegraphics[width=5.9in]{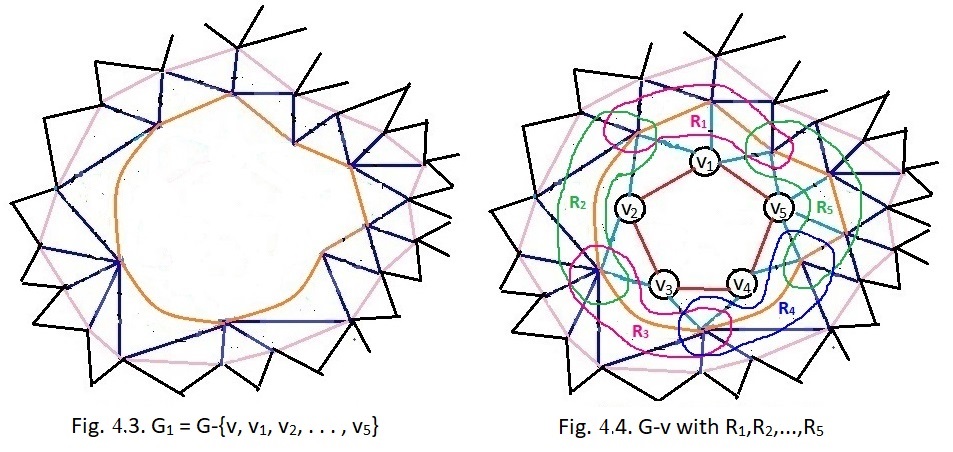}}
 \end{figure}
 \begin{figure}[ht]
 	\centerline{\includegraphics[width=3.2in]{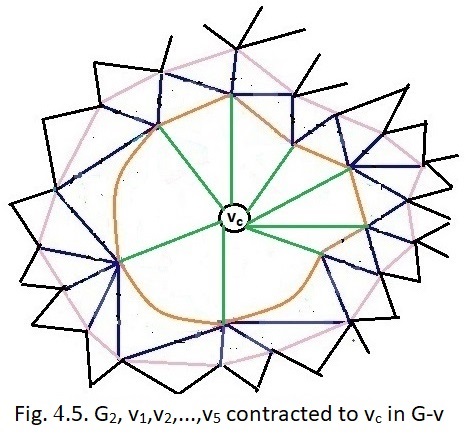}}
 \end{figure}
 
\item[{\rm 3.}] In {\color{red} any proper coloring of} $G-v$, color of $v_i$ should be different from the elements of $R_i$ and colors of $v_{i-1}$ and $v_{i+1}$. See Figures 4.4.

\item[{\rm 4.}] Two adjacent $R_i$s, say $R_i$ and $R_{i+1}$, have at least one common element corresponding to the color of the vertex which is adjacent to both $v_i$ and $v_{i+1}$ in $G-v$ since $G$ is a maximal planar graph. See Figure 5.1.

\item[{\rm 5.}] In $G-v$, there is only one vertex which is adjacent to both $v_i$ and $v_{i+1}$. Otherwise it comes under case-1. See Figure 5.2.
\begin{figure}[ht]
	\centerline{\includegraphics[width=4.5in]{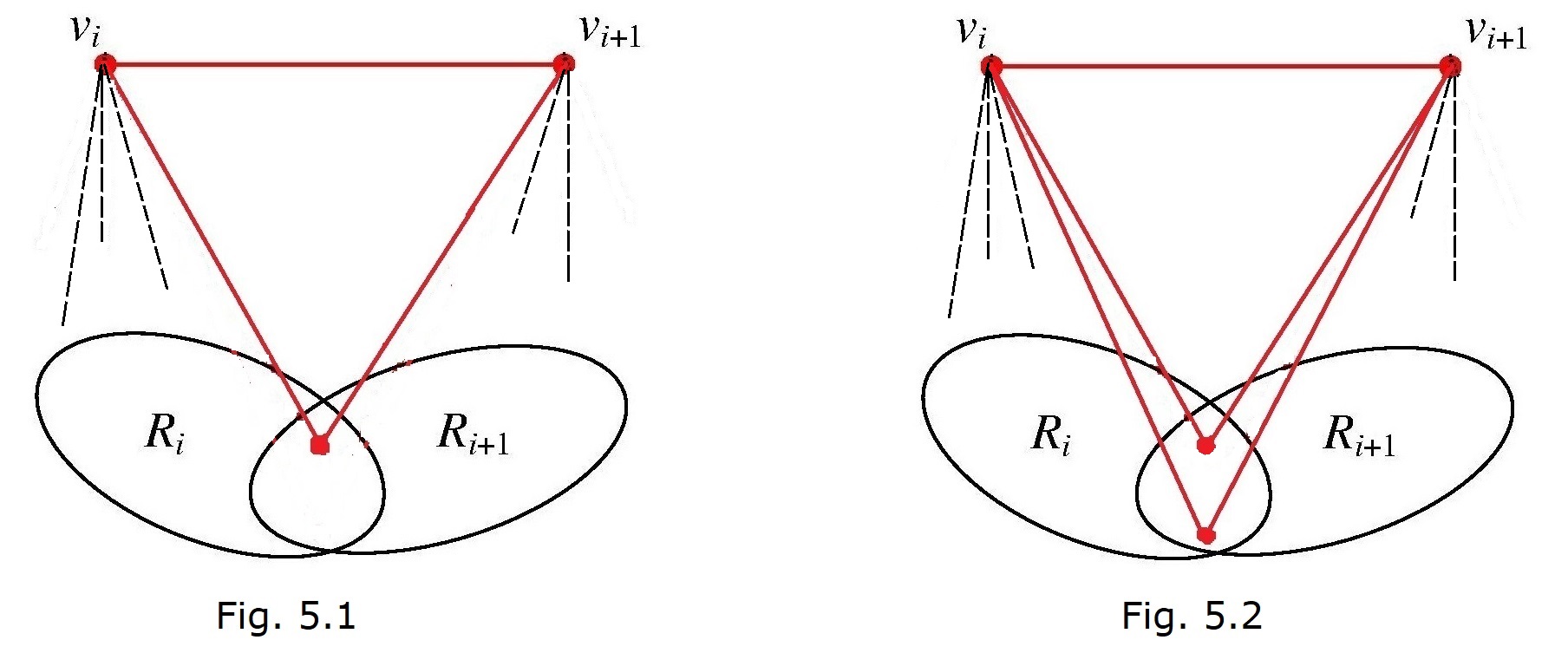}}
\end{figure}

\item[{\rm 6.}] Two adjacent $R_i$s, say $R_i$ and $R_{i+1}$, each having all the three colors {\color{blue} of $R_1$ $\cup$ $R_2$ $\cup$ $\ldots$ $\cup$ $R_{d(v)}$ leads to $v_i$ and $v_{i+1}$ taking $4^{th}$ and $5^{th}$ colors in $G-v$ and so we should avoid such colorings.} See Figure 5.1.

\item[{\rm 7.}] In each figure number just out side $R_i$ corresponds to the color of vertex $v_i$ in $G-v,$ {\color{blue} $1 \leq i \leq d(v)$.}

\item[{\rm 8.}]  While considering different possible colorings of graph $G-v$, we avoid repetition of equivalent configurations. 

Here, \textit{two configurations or figures (of colorings) $F_1$ and $F_2$} (of a graph $H_1$) are called \textit{equivalent} if mutual interchange of numbers (colors) 1, 2, 3, 4 in one figure give rise to the other figure. That is there exists a permutation between colors of color classes of $F_1$ and $F_2$. In this case, we denote it by $F_1$ $\cong$ $F_2$. See Figures 6.1 to 6.3.

\begin{figure}[ht]
	\centerline{\includegraphics[width=4.8in]{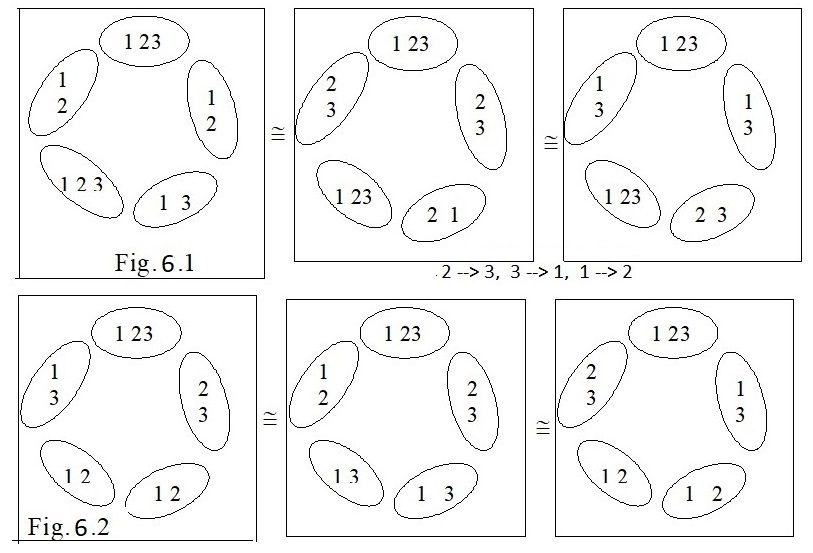}}
	\centerline{\includegraphics[width=4.9in]{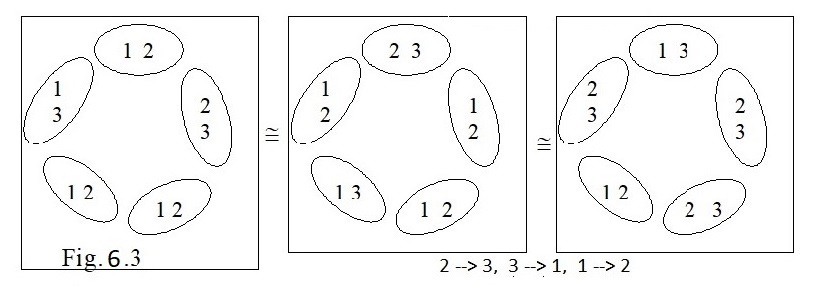}}
\end{figure}
\end{enumerate}

By all the above conditions, when $d(v)$ $=$ $4,$ it is always possible to find {\color{blue} proper} colorings of $G-v$ in which $G-v$ is 4 colorable and $v_1$, $v_2$, $v_3$ and $v_4$ take at the most three colors out of the 4 colors. See Figures 7.1 to 7.5. They are the {\color{red} proper} colorings with more colors under $d(v)$ $=$ $4$ (in the sense, if the number of elements of $R_is$ are reduced further, then it is easy to see that $v_1$, $v_2$, $v_3$ and $v_4$ require at the most 3 colors out of the 4 colors) {\color{red} under the above conditions}. {\color{red} Thus, in the case of $d(v)$ = 4}, by applying a color out of the 4 colors and different from that of $v_1$, $v_2$, $v_3$ and $v_4$ to $v$, the graph $G$ is 4 colorable.

\begin{enumerate}
\item[{\rm 9.}] When {\color{blue}{\bf $d(v) = 5$}}, $d(v_i) \geq d(v)$ = 5 and each $R_i$ contains at least two distinct elements, $1 \leq i \leq 5$.  Otherwise, {\color{red} we get $d(v_i) \leq 4$ in $G,$} $1 \leq i \leq 5.$ See Figures 8.1 and 8.2.
\begin{figure}[ht]
	\centerline{\includegraphics[width=4.3in]{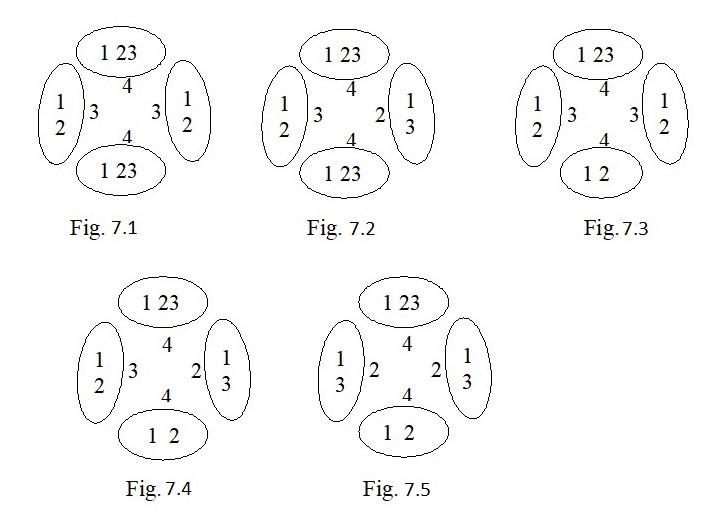}}
	\centerline{\includegraphics[width=4.5in]{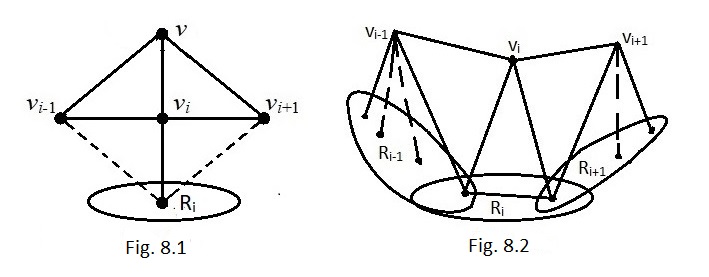}}
\end{figure}
\item[{\rm 10.}] Using the above conditions, we obtain Figures 9.1 to 9.88 when $d(v)$ = 5. See  Figures 9.1 to 9.88. The number of these figures can be reduced to 23 using condition-8 and only for clarity 88 figures (type of colorings) are considered. 
\begin{figure}[ht]
	\centerline{\includegraphics[width=4.3in]{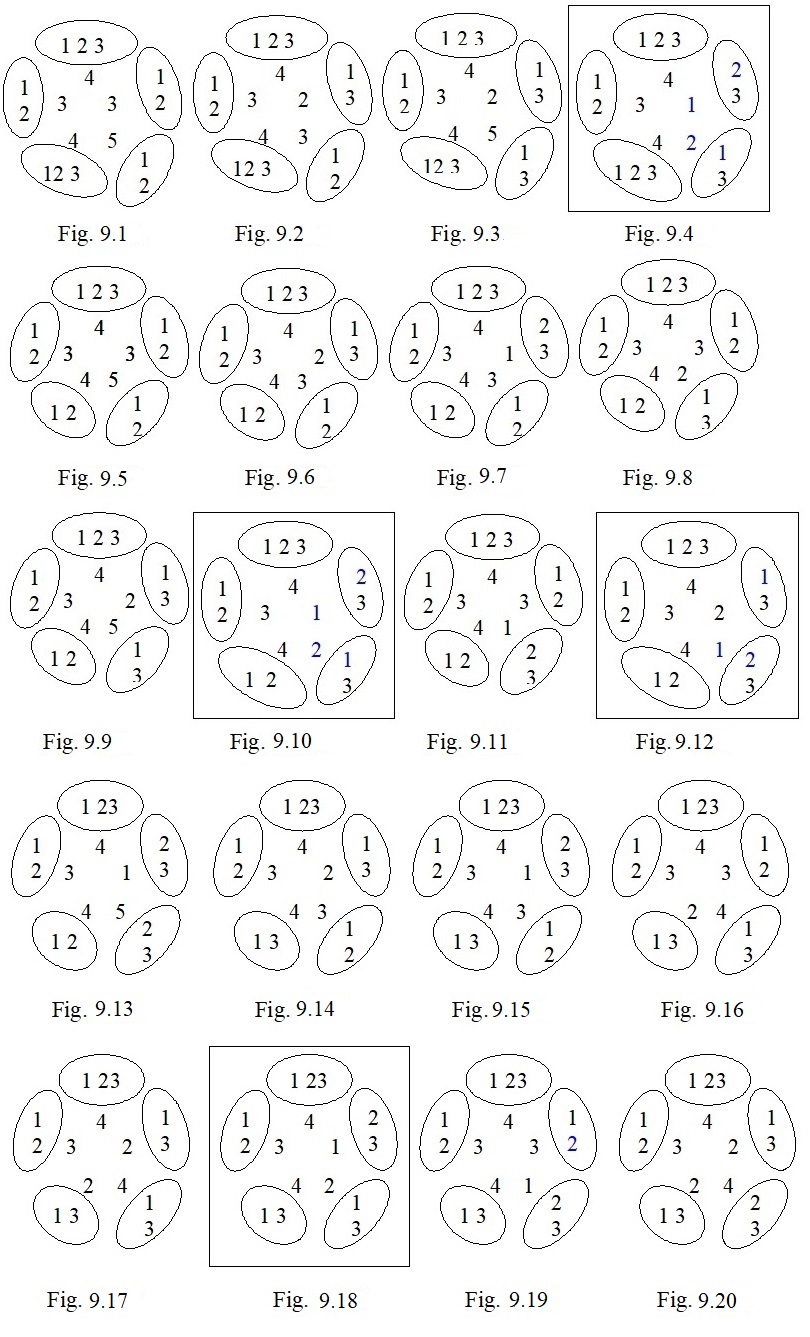}}
\end{figure}
\begin{figure}[ht]
	\centerline{\includegraphics[width=4.3in]{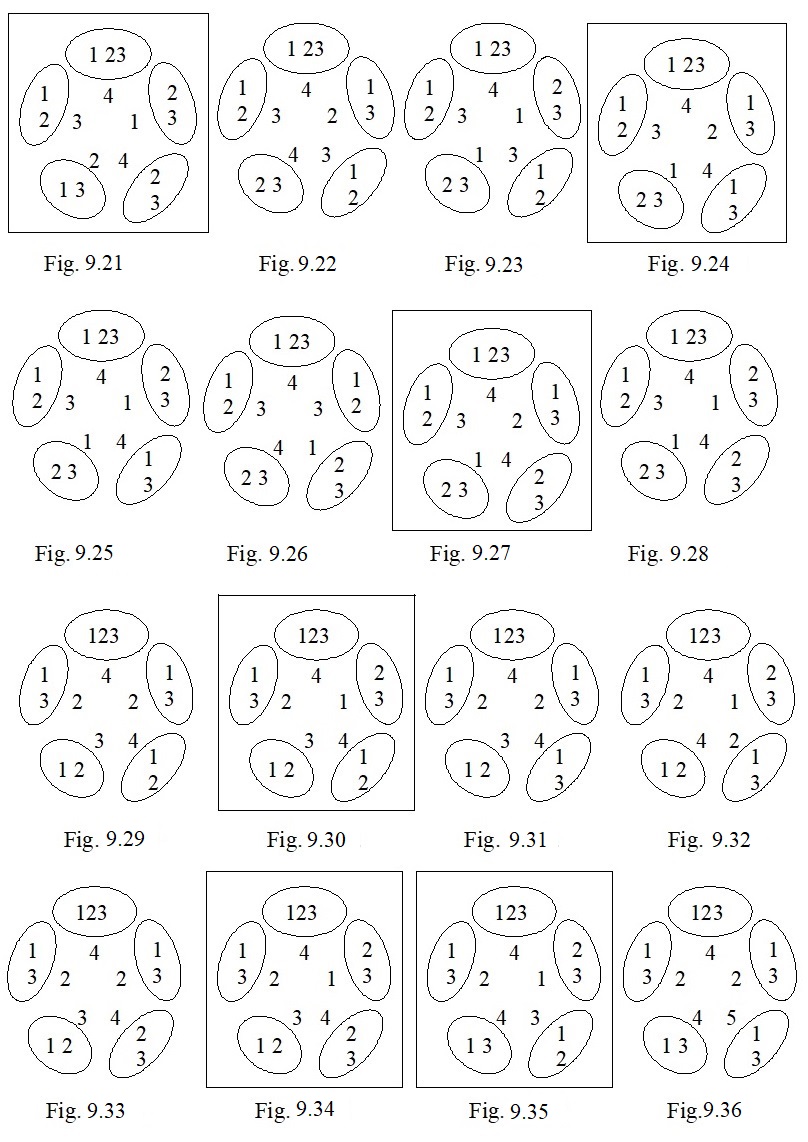}}
\end{figure}
\begin{figure}[ht]
	\centerline{\includegraphics[width=4.3in]{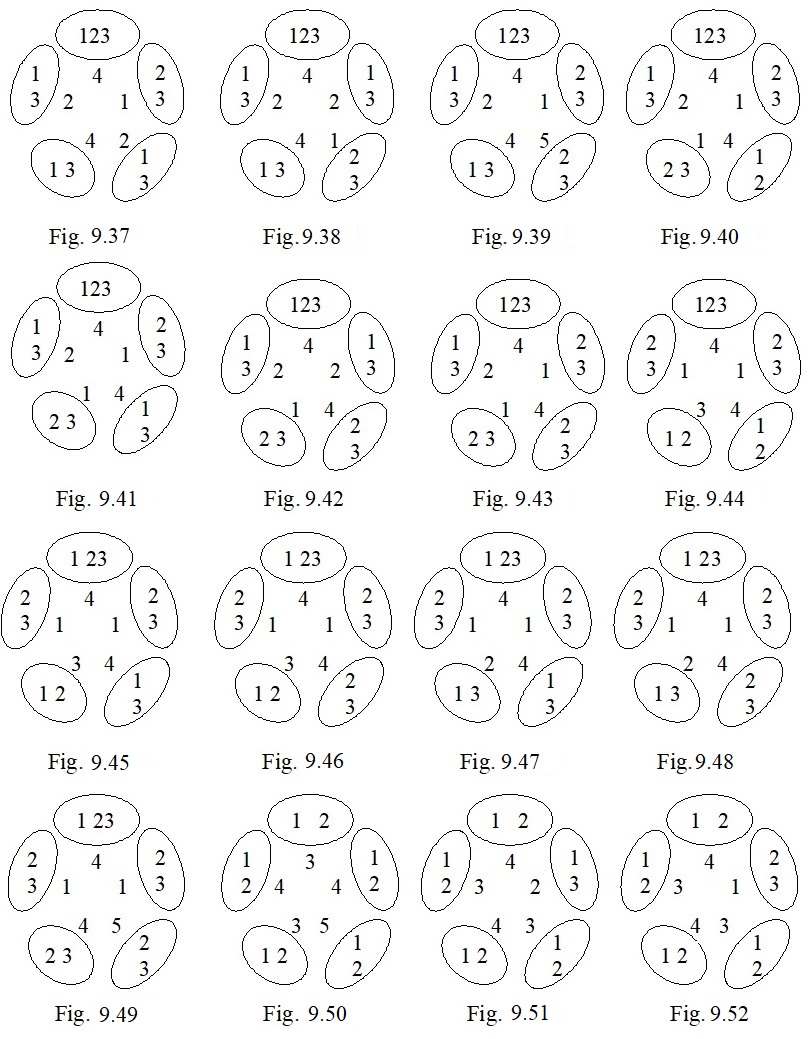}}
\end{figure}
\begin{figure}[ht]
	\centerline{\includegraphics[width=4.3in]{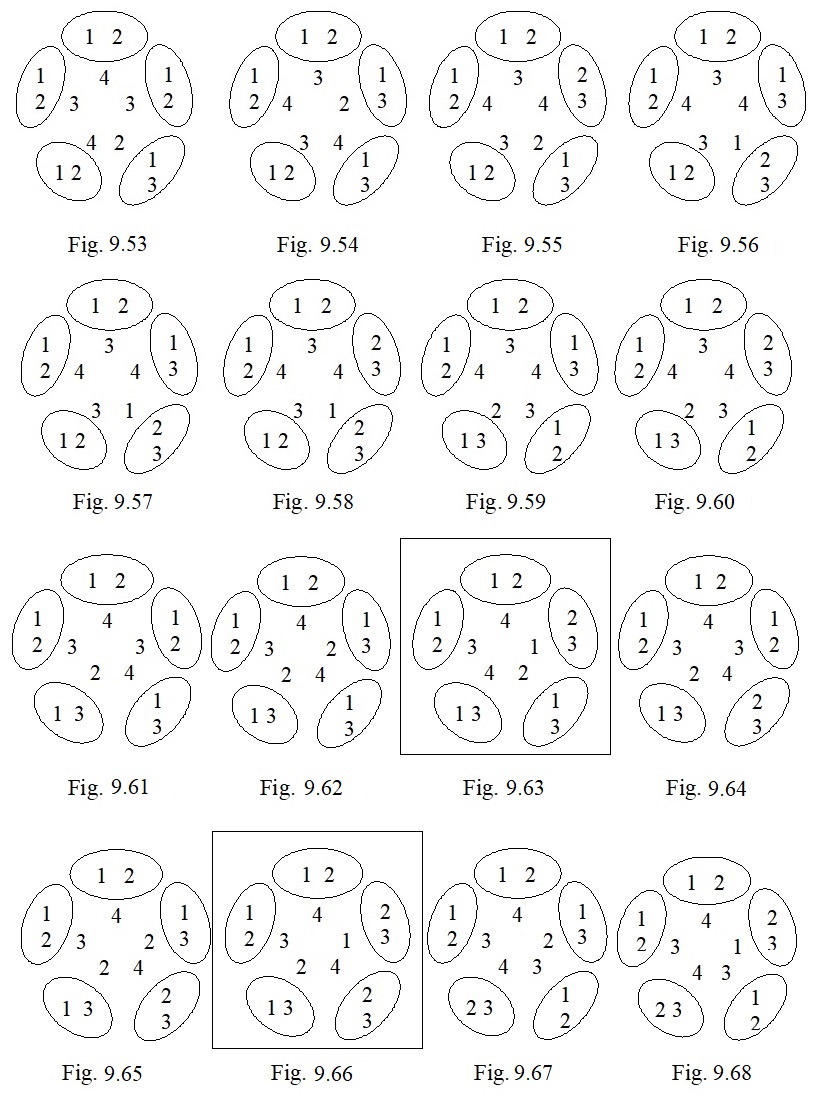}}
\end{figure}
\begin{figure}[ht]
	\centerline{\includegraphics[width=4.5in]{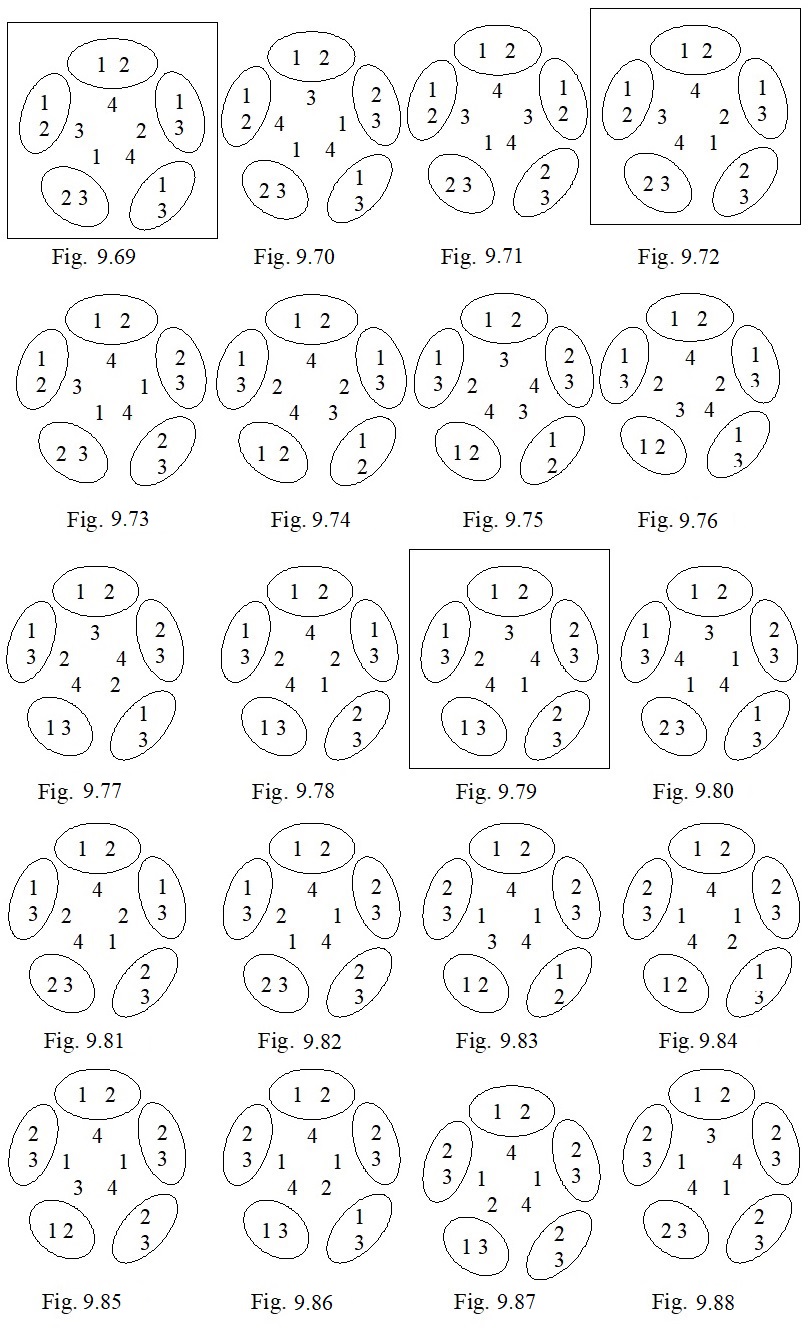}}
\end{figure}
The different and equivalent figures (colorings) under Figure 9 are

\begin{enumerate}
\item [{\rm (1)}] Fig. 9.1; 

\item [{\rm (2)}] Fig. 9.2; 

\item [{\rm (3)}] Fig. 9.3; 

\item [{\rm (4)}] Fig. 9.4; 

\item [{\rm (5)}] Fig. 9.5 $\cong$ Fig. 9.36 $\cong$ Fig. 9.49; 

\item [{\rm (6)}] Fig. 9.6 $\cong$ Fig. 9.7 $\cong$ Fig. 9.17 $\cong$ Fig. 9.28 
$\cong$ Fig. 9.37 $\cong$  Fig. 9.43; 

\item [{\rm (7)}] Fig. 9.8 $\cong$ Fig. 9.11 $\cong$ Fig. 9.31 $\cong$ Fig. 9.38 $\cong$ Fig. 9.46 $\cong$ Fig. 9.48; 

\item [{\rm (8)}] Fig. 9.9 $\cong$ Fig. 9.13 $\cong$ Fig. 9.39; 

\item [{\rm (9)}] Fig. 9.10 $\cong$ Fig. 9.12 $\cong$ Fig.9.21 $\cong$ Fig. 9.24 
$\cong$ Fig. 9.34 $\cong$  Fig. 9.35; 

\item [{\rm (10)}] Fig. 9.14 $\cong$ Fig. 9.23 $\cong$ Fig. 9.41; 

\item [{\rm (11)}] Fig. 9.15 $\cong$ Fig. 9.20 $\cong$ Fig. 9.22 $\cong$ Fig. 9.25 
$\cong$ Fig. 9.32 $\cong$  Fig. 9.40; 

\item [{\rm (12)}] Fig. 9.16 $\cong$ Fig. 9.26 $\cong$ Fig. 9.29 $\cong$ Fig. 9.42 $\cong$ Fig. 9.44 $\cong$  Fig. 9.47;

\item [{\rm (13)}] Fig. 9.18 $\cong$ Fig. 9.27 $\cong$ Fig. 9.30; 

\item [{\rm (14)}] Fig. 9.19 $\cong$ Fig. 9.33 $\cong$ Fig. 9.45; 

\item [{\rm (15)}] Fig. 9.50; 

\item [{\rm (16)}] Fig. 9.51 $\cong$ Fig. 9.52 $\cong$ Fig. 9.53 $\cong$ Fig. 9.88; 

\item [{\rm (17)}] Fig. 9.54 $\cong$ Fig. 9.58 $\cong$ Fig. 9.61 $\cong$ Fig. 9.62  $\cong$ Fig. 9.71 $\cong$  Fig. 9.73; 

\item [{\rm (18)}] Fig. 9.55 $\cong$ Fig. 9.56 $\cong$ Fig. 9.57 $\cong$ Fig. 9.64  $\cong$ Fig. 9.77 $\cong$  Fig. 9.82; 

\item [{\rm (19)}] Fig. 9.59 $\cong$ Fig. 9.68 $\cong$ Fig. 9.74 $\cong$ Fig. 9.76  $\cong$ Fig. 9.83 $\cong$  Fig. 9.85; 

\item [{\rm (20)}] Fig. 9.60 $\cong$ Fig. 9.67 $\cong$ Fig. 9.75 $\cong$ Fig. 9.78  $\cong$ Fig. 9.87; 

\item [{\rm (21)}] Fig. 9.63 $\cong$ Fig. 9.66 $\cong$ Fig. 9.69 $\cong$ Fig. 9.72  $\cong$ Fig. 9.79; 

\item [{\rm (22)}] Fig. 9.65 $\cong$ Fig. 9.70 $\cong$ Fig. 9.81 $\cong$ Fig. 9.86 and 

\item [{\rm (23)}] Fig. 9.80 $\cong$ Fig. 9.84. 
\end{enumerate}

\begin{enumerate}
\item[{\rm 11.}] Colorings corresponding to Figures 9.1, 9.3, 9.5, 9.9, 9.13, 9.36, 9.39, 9.49 and 9.50 are avoidable colorings of $G-v$ since each one contains $5^{th}$ color. And all other figures (type of colorings) under Figure 9, other than figures each inside a square, correspond to possible {\color{blue} proper} colorings with which $G$ is colorable with $\leq$ 4 colors.

\item[{\rm 12.}] Figure inside a square, under Figure 9, corresponds to {\color{red} proper coloring (of $G-v$ with four colors)} in which $v_1$, $v_2$, $v_3$, $v_4$, $v_5$ require colors 1, 2, 3 and 4. These figures are 9.4, 9.10, 9.12, 9.18, 9.21, 9.24, 9.27, 9.30, 9.34, 9.35, 9.63, 9.66, 9.69, 9.72 and 9.79 and using equivalent colorings these comes under four different figures, namely,

\begin{enumerate}
	\item [{\rm (a)}] Fig. 10.1 = Fig. 9.4;

\item [{\rm (b)}] Fig. 10.2 = Fig. 9.10 $\cong$ Fig. 9.12 $\cong$ Fig. 9.21 

~ \hfill $\cong$ Fig. 9.24 $\cong$ Fig. 9.34 $\cong$ Fig. 9.35;

\item [{\rm (c)}] Fig. 10.3 = Fig. 9.18 $\cong$ Fig. 9.27 $\cong$ Fig. 9.30 and

\item [{\rm (d)}] Fig. 10.4 = Fig. 9.63 $\cong$ Fig. 9.66 $\cong$ Fig. 9.69 

~ \hfill $\cong$ Fig. 9.72 $\cong$ Fig. 9.79.
\end{enumerate}

We consider them separately and by applying contraction we find out alternate colorings to these. Figures 10.1 to 10.4 are also given in Figure 11. {\color{red} See Figures 10.1 to 10.4 and 11.}
\begin{figure}[ht]
	\centerline{\includegraphics[width=4.8in]{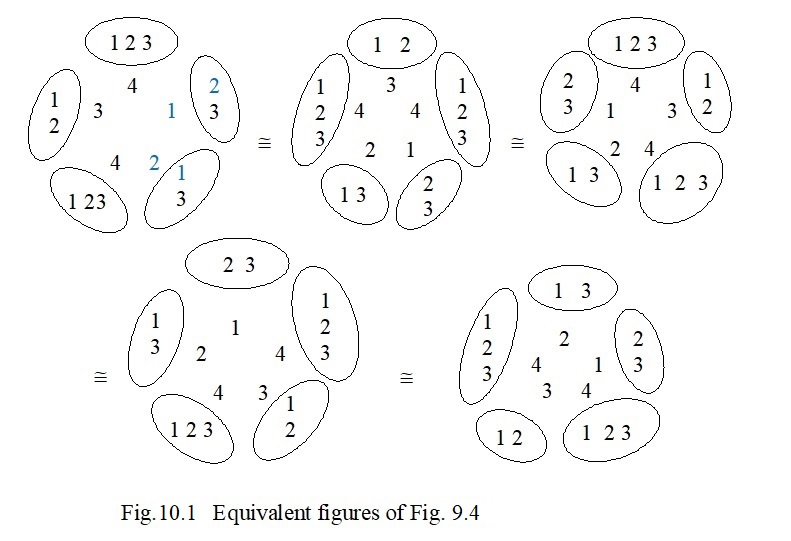}}
	\centerline{\includegraphics[width=4.8in]{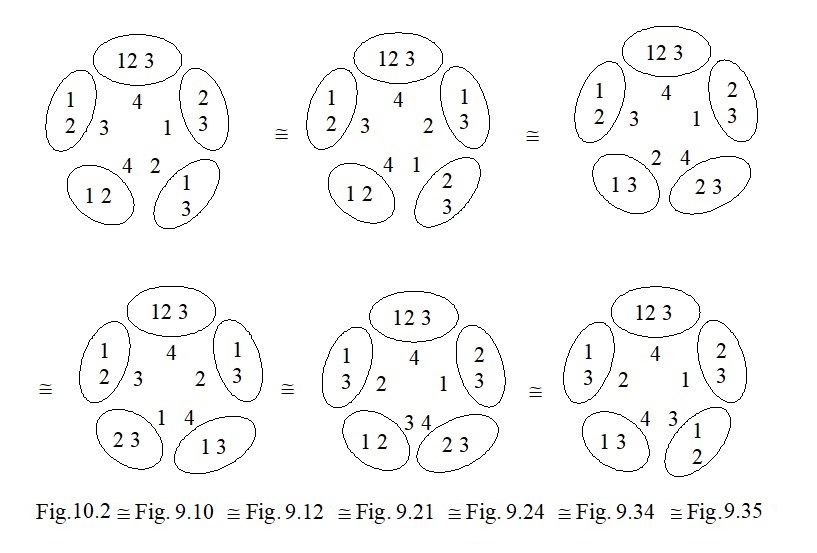}}
\end{figure}
\begin{figure}[ht]
	\centerline{\includegraphics[width=4.8in]{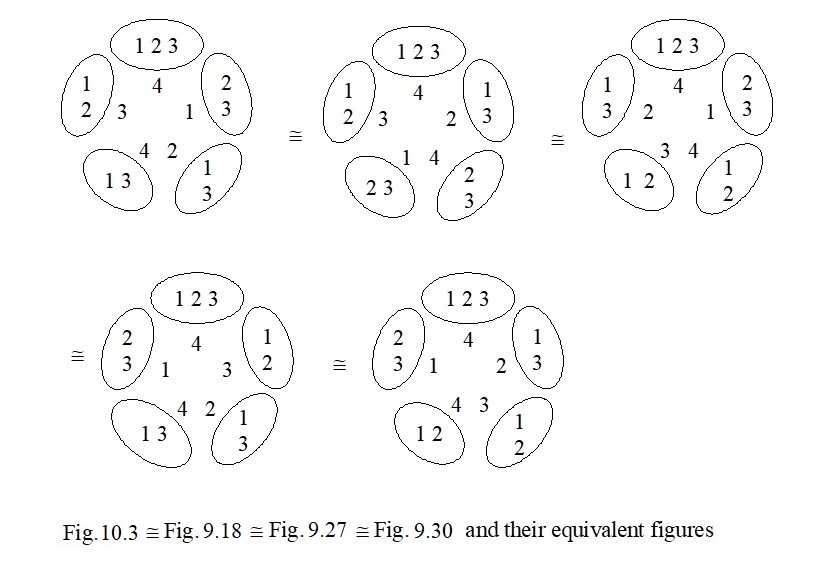}}
	\centerline{\includegraphics[width=4.8in]{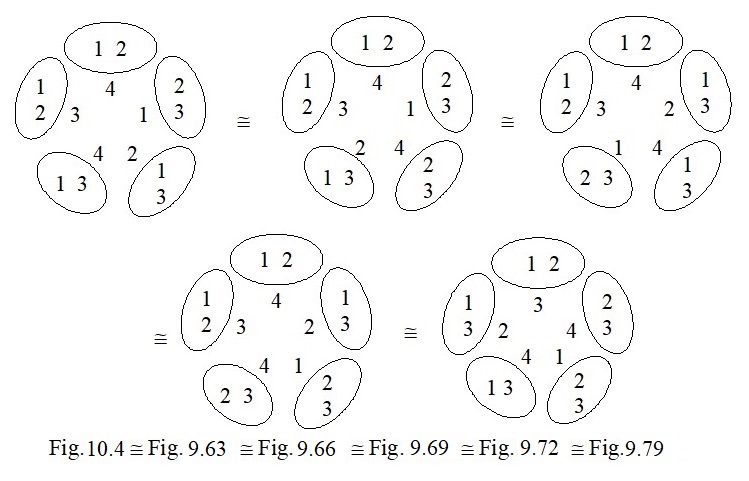}}
\end{figure}
\item[{\rm 13.}] Each of Figures 10.1 to 10.4 contains adjacent vertices $v_i$ and $v_{i+1}$ of the cycle $(v_1 ~v_2~ v_3~ v_4 ~v_5)$  {\color{red}  such that ($R_i$ and $R_{i+1}$ are adjacent in $G-v$ and)} $R_{i} \neq R_{i-1}$ and $R_{i+1} \neq R_{i+2}$, {\color{blue} $1 \leq i-1,i,i+1,i+2 \leq 5$.} And corresponding $R_i$ and $R_{i+1}$ in Figures 10.1 to 10.4 are as follows.

\noindent
In Figure 10.1.: {\color{blue} Pair of such adjacent $R_i$s in the figure}, $1 \leq i \leq 5$.
\begin{enumerate}
\item $R_1$ = $\{1, 2, 3\}$ and $R_2$ = $\{1, 2\}$ (Here, $R_1 \neq R_5$ and $R_2 \neq R_3$.);
\item $R_2$ = $\{1, 2\}$ and $R_3$ = $\{1, 2, 3\}$ (Here, $R_2 \neq R_1$ and $R_3 \neq R_4$.);
\item $R_3$ = $\{1, 2, 3\}$ and $R_4$ = $\{1, 3\}$ (Here, $R_3 \neq R_2$ and $R_4 \neq R_5$.);
\item $R_4$ = $\{1, 3\}$ and $R_5$ = $\{2, 3\}$ (Here, $R_4 \neq R_3$ and $R_5 \neq R_1$.);
\item $R_5$ = $\{2, 3\}$ and $R_1$ = $\{1, 2, 3\}$ (Here, $R_5 \neq R_4$ and $R_1 \neq R_2$.).
\end{enumerate}
\noindent
In Figure 10.2.:
\begin{enumerate}
\item $R_2$ = $\{1, 2\}$ and $R_3$ = $\{1, 2\}$ (Here, $R_2 \neq R_1$ and $R_3 \neq R_4$);
\item $R_4$ = $\{1, 3\}$ and $R_5$ = $\{2, 3\}$ (Here, $R_4 \neq R_3$ and $R_5 \neq R_1$);
\item $R_5$ = $\{2, 3\}$ and $R_1$ = $\{1, 2, 3\}$ (Here, $R_5 \neq R_4$ and $R_1 \neq R_2$.).
\end{enumerate}
\noindent
In Figure 10.3.:
\begin{enumerate}
\item $R_1$ = $\{1, 2, 3\}$ and $R_2$ = $\{1, 2\}$ (Here, $R_1 \neq R_5$ and $R_2 \neq R_3$.);
\item $R_3$ = $\{1, 3\}$ and $R_4$ = $\{1, 3\}$ (Here, $R_3 \neq R_2$ and $R_4 \neq R_5$.);
\item $R_5$ = $\{2, 3\}$ and $R_1$ = $\{1, 2, 3\}$ (Here, $R_5 \neq R_4$ and $R_1 \neq R_2$.).
\end{enumerate}
\noindent
In Figure 10.4.:
\begin{enumerate}
\item $R_1$ = $\{1, 2\}$ and $R_2$ = $\{1, 2\}$ (Here, $R_1 \neq R_5$ and $R_2 \neq R_3$.);
\item $R_3$ = $\{1, 3\}$ and $R_4$ = $\{1, 3\}$ (Here, $R_3 \neq R_2$ and $R_4 \neq R_5$.). 

See Figure 11.
\end{enumerate}

\item[{\rm 14.}]  In each of Figures 10.1 to 10.4 of $G-v$, make contraction on one set of such edges $v_{i-1} v_i$ and $v_{i+1} v_{i+2}$, {\color{blue} $1 \leq i-1,i,i+1,i+2 \leq 5$.} Let the corresponding contracted graph of $G-v$ be $G_3$, contracted vertices of the contracted edges be $v_{i-1, i}$ and $v_{i+1, i+2}$ and the contracted $R_j$s be $R_{i-1, i}$ and $R_{i+1, i+2}$, respectively, {\color{blue} $1 \leq i-1,i,i+1,i+2,j \leq 5$.}   

Clearly, the cycle $(v_1 ~v_2~ v_3 ~v_4~ v_5)$ is contracted to a cycle of length three, two of its vertices being $v_{i-1, i}$ and $v_{i+1, i+2}$ in $G_3$ and $R_{i-1, i}$ = $R_{i-1} \cup R_{i}$ =  $\{1, 2, 3\}$ = $R_{i+1} \cup R_{i+2}$ = $R_{i+1, i+2}$. For an illustration, see Figures 12.1 and 12.2. 
\begin{figure}[ht]
	\centerline{\includegraphics[width=3.8in]{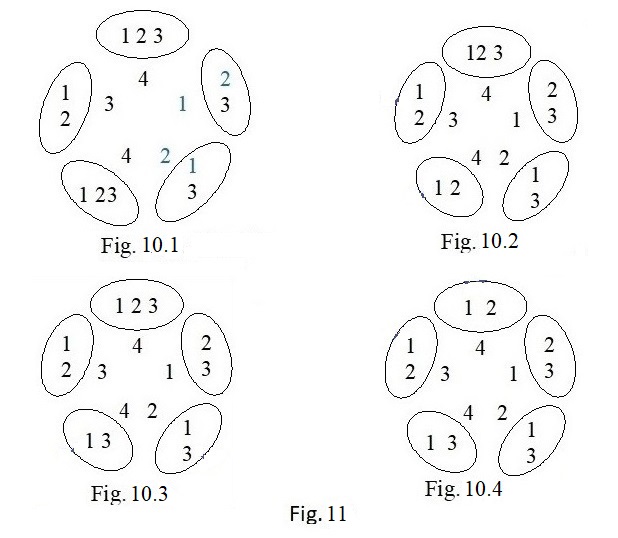}}
	\centerline{\includegraphics[width=4.3in]{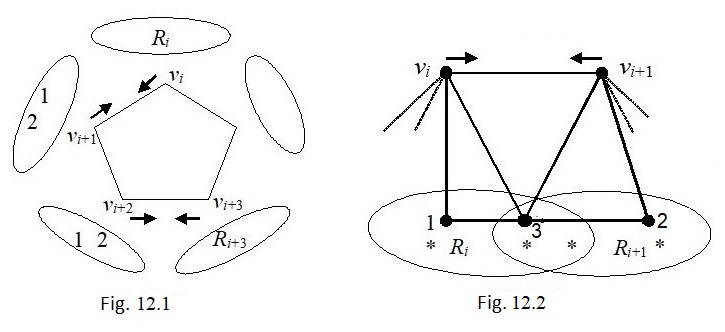}}
\end{figure}
Hence, $v_{i-1, i}$ and $v_{i+1, i+2}$ require $4^{th}$ and $5^{th}$ colors in $G_3$ and so the coloring is {\color{blue} an avoidable}   coloring to $G_3$ {\color{blue} whose order is $< n$}. By assumption, there exists an alternate coloring to $G_3$ (and thereby to $R_1$, $R_2$, $R_3$, $R_4$, $R_5$ and thereby to $G_1$ (of $G_3$ as well as of $G-v$) and thereby to $G-v$) in which $R_i$s should not satisfy condition-13 and so it should be different from Figures 10.1 to 10.4. 

In Figure 13, red lines with arrow marks indicate contraction made on one set of such edges $v_{i-1} v_i$ and $v_{i+1} v_{i+2}$ of Figures 10.1 to 10.4. See Figure 13.
\begin{figure}[ht]
	\centerline{\includegraphics[width=3.5in]{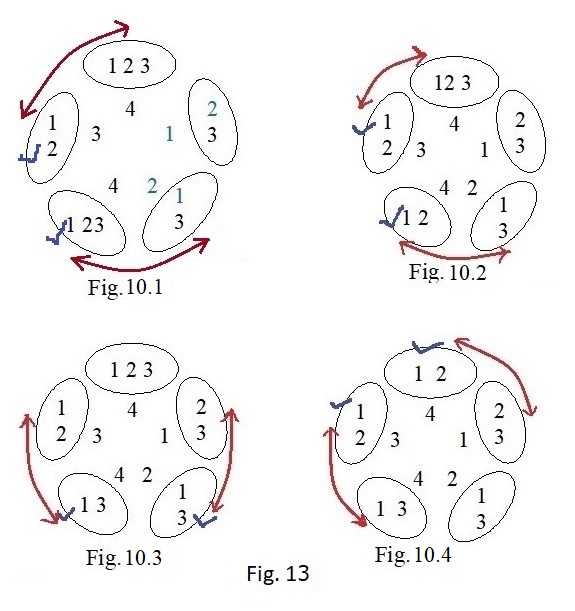}}
\end{figure}

Now, the question is whether such an alternate coloring exists without satisfying condition-13. Yes. Coloring corresponding to Fig. 9.51 $\cong$ Fig. 9.52 $\cong$ Fig. 9.53 $\cong$ Fig. 9.88 is an alternate coloring.

This implies, among all possible {\color{blue} proper colorings of $G-v$ with 4 colors,} colorings in which the cycle $(v_1 ~v_2 ~v_3 ~v_4 ~v_5)$ taking 4 colors are having alternate colorings other than the colorings corresponding to Figures 10.1 to 10.4.  

Thus, there exist {\color{blue} proper} colorings to $G-v$ such that $G-v$ is 4 colorable while its vertices $v_1$, $v_2$, $v_3$, $v_4$ and $v_5$ take less than 4 colors out of the 4 colors. By assigning a color out of the four colors, other than the colors of $v_1$, $v_2$, $v_3$, $v_4$ and $v_5$ to $v$, the graph $G$ is 4 colorable. Thus, when $d(v)$ = 5, the graph $G$ of order $n+1$ is four colorable.
\end{enumerate}
\end{enumerate}

Hence, the claim is true for $G$, maximal planar graph of order $n+1$. And thereby the result is true for $H$, simple planar graph which is a subgraph of $G$ and of order $n+1$. Thus, we could prove the result for $n+1$. Therefore, by mathematical induction, the theorem is true for all values of $n$, $n\in\mathbb{N}$. That is any simple planar graph of order $n$ is four colorable, $n\in\mathbb{N}$.
\end{proof}

\vspace{.2cm}
\noindent
{\bf Acknowledgement}\quad I express my deep sense of gratitude to Prof. Brian Alspach, University of Newcastle, Australia, who helped wholeheartedly at every stage to bring the manuscript to this stage. I  express my sincere thanks to Prof. S. Krishnan (late), Prof. V. Mohan and Prof. R. Aravamudhan (late), Thiagarajar College of Engineering, Madurai, India; Prof. Lowell W Beineke, Purdue University, U.S.A.; Dr. K. Vareethaiah and Dr. S. Amirthaiyan, St. Jude's College, Thoothoor, India; Dr. Oscar Fredy, Royal Liverpool University Hospital, Liverpool, U.K. and Dr. M. I. Jinnah, Mr. Albert Joseph and Dr. L. John, University of Kerala, Trivandrum, Kerala, India for their help and guidance. I also express my sincere gratitude to unknown referees, who helped at different stages  and to the Central University of Kerala,  Kasaragod, India - 671 320 and St. Jude's College, Thoothoor, TN, India - 629 176 for providing facilities.

\end{document}